\documentclass{article}
\usepackage{amsmath,amsthm,amsfonts,amscd,amssymb,eucal,latexsym}

\begin{document}

\title{The representation category of the Woronowicz
quantum group $S_\mu U(d)$ as a braided tensor $C^*$--category}
\author{Claudia Pinzari\\
Dipartimento di Matematica, Universit\`a di Roma La Sapienza\\
00185--Roma, Italy}
\date{}
\maketitle

\begin{abstract} An abstract characterization of the representation
category of the Woronowicz twisted $SU(d)$ group is given, generalizing 
analogous results known in the classical case.
\end{abstract}

\begin{section}{Introduction} 

In \cite{DR} Doplicher and Roberts proved a duality theory 
for compact groups which allows one to recognize the representation category of a group
$G$  among all tensor  $C^*$--categories as those endowed 
with a symmetry of the permutation group, and for which each object has a conjugate.
\smallskip

This paper is part of  the 
program of generalizing the Doplicher--Roberts  duality theory to
Woronowicz 
compact quantum groups. As Woronowicz showed in \cite{Wo}, the
representation categories in question
are precisely the
 tensor $^*$--subcategories of Hilbert space categories for which each object has a conjugate
in the sense of \cite{LR}. Therefore, at this level of  generality, 
we lose any symmetry in the category.\smallskip

However, motivated by Quantum Field Theory in low dimensions, and with the final aim
of a possible application to the physical situation, we are  interested in 
studying tensor $C^*$--categories which have  a braided symmetry: a
representation of the braid group
${\Bbb B}_n$ in the intertwiner spaces  $(\rho^n,\rho^n)$ between 
 powers of each object $\rho$
of the category which preserves the tensor structures. \smallskip

Here, for example, we look at the
Woronowicz quantum group
$S_\mu U(d)$. Its representation category
contains, among its intertwiners, a representation of the braid group: the
Jimbo--Woronowicz representation \cite{Ji}, \cite{Wo}. This is of a very
particular
form, as the basic intertwiner $g_q$ is selfadjoint  with
only two eigenvalues: $\mu^2$ and $-1$.
This means 
that the braided symmetry in fact
factors through a representation of the Hecke algebra $H_\infty(\mu^2)$.
\smallskip

We show that actually this Hecke symmetry satisfies more: its basic
intertwiner $g_q$
can be suitably normalized so that the resulting braided symmetry $\varepsilon$
makes  the representation category of $S_\mu U(d)$ into a {\sl braided tensor category}, 
in the sense that 
$$\varepsilon(g_1\dots g_m)T\otimes 1_H=1_H\otimes T\varepsilon(g_1\dots g_n)$$
whenever $T\in (H^{\otimes n}, H^{\otimes m})$ is an intertwiner
 in the representation category of
$S_\mu U(d)$ between tensor powers of the defining representation of $S_\mu U(d)$ on the 
$d$--dimensional Hilbert space $H$ (Cor. 5.2).\smallskip

Our interest in braiding is  motivated by the fact 
that the tensor $C^*$--categories arising from low dimensional 
QFT are braided, albeit with a {\it unitary} braiding.
\smallskip

The aim of this paper is to characterize the representation category of $S_\mu U(d)$
among all braided tensor $C^*$--categories with conjugates (Theorem
6.2).\smallskip

In section 2 we review the notion of Hecke algebra and the properties we
shall need.
In particular, we introduce special cases of the Young 
symmetrizers and antisymmetrizers
which have appeared  in the literature 
concerning representations of Hecke algebras (see, e.g., \cite{DJ}).
\smallskip

In section 3 we study Hecke symmetries in tensor $C^*$--categories and we
generalize
results previously known for permutation symmetries due to Doplicher
and Roberts \cite{DRduals}. The main result of 
this section is the generalization of the theorem about the restriction of
the `statistics parameter' values to the values
$\lambda_d:=q^d\frac{q-1}{q^d-1}$, with $d\in{\Bbb Z}$. Here we consider
both the case where
$q$ is real and the case where $q$ is a root of unity. 
\smallskip

We show that each
of these allowed values determines uniquely the kernel
of the 
Hecke algebra  representations given by the symmetry: If $q>0$ and $d$ 
is negative, then the
kernel is the same as the kernel of the
Jimbo--Woronowicz
representation on $H\otimes H$ with $H$ a $d$--dimensional Hilbert space.
Symmetries where this situation occurs will be then called of dimension
$d$.
\smallskip

In the root of unity case we find
out that
the given Hecke algebra representation with parameter $\lambda_d$ has the
same kernel as the 
Wenzl's representation $\pi^{(d,m)}$
(see Theorem 3.3 for a precise statement).\smallskip

Later on we specialize to the case where $q>0$. 
In section 5 we introduce the basic intertwiner, namely the quantum
determinant,
$S$ of $S_\mu U(d)$, 
with $\mu=\sqrt{q}$, and we
compute the conjugate representation of the defining
 representation (Theorem 5.5). For future reference,
we perform the computations in a slight more generality, considering also
quantum determinants  of proper subspaces of $H$. In doing so we
discover that
in the case where $q\neq1$ there exist more  
left inverses of $H$ which are faithful
on the image of the Jimbo--Woronowicz representation than in the classical
case (Cor.
5.2). We realize also that these smaller quantum determinants  do
not
satisfy the afore mentioned braid relation.\smallskip

In section 6 we prove the main result which characterizes the
representation category of $S_\mu U(d)$ among tensor $C^*$--categories by
means of its Hecke symmetry
and the conjugate representation of the fundamental representation.

\end{section}

\begin{section}{Preliminaries on the Hecke algebra}

For any integer $n$ 
let 
$B_n$ denote the braid group with generators
$g_i$, $i=1,\dots, n-1$ and relations
$$g_ig_j=g_jg_i\quad i,j: |i-j|>1\eqno(2.1)$$
$$g_ig_{i+1}g_i=g_{i+1}g_ig_{i+1}.\eqno(2.2)$$
 For any  complex number $q$, let $H_n(q)$ denote the 
Hecke algebra of order $n$, namely the quotient of the complex group
algebra
${\Bbb C}[{\Bbb B}_n]$ by the relations
$${g_i}^2=(q-1)g_i+q, i=1,\dots,n-1.\eqno(2.3)$$

We also consider the inductive limit $H_\infty(q)$ of $(H_n(q),
\iota_{n,m})$,
where, for $m>n$, $\iota_{n,m}$ takes the element $g_i\in H_n(q)$ to 
the element $g_i\in H_m(q)$.
Therefore the map $g_i\to g_{i+1}$ extends uniquely to a monomorphism
$\sigma$ of
$H_\infty(q)$ into itself.

Let us introduce, for $q\neq -1$, the elements $e_i=\frac{1+g_i}{q+1}$.
Then  $(2.3)$ is equivalent to the fact that $e_i$ is an idempotent:
${e_i}^2=e_i$. Therefore $g_i$ is a linear combination of two mutually 
orthogonal
idempotents:
$$g_i=qe_i-(1-e_i).$$

We define, for $q\neq0$, special   elements $A_k$, $k\in{\Bbb
Z}-\{0\}$, of
$H_\infty(q)$:
$$A_1=A_{-1}=1,$$
$$A_{n+1}=\sigma(A_n)+g_1\sigma(A_n)+\dots+g_{n}\dots g_1\sigma(A_n),$$
$$A_{-n-1}=\sigma(A_{-n})-q^{-1} g_1\sigma(A_n)+\dots+
(-q^{-1})^{n}g_{n}\dots
g_1\sigma(A_{-n}).$$

Notice that both $A_2=(1+q)e_1$ and $A_{-2}=(1+q^{-1})(1-e_1)$.
These elements are special cases of symmetrization and antisymmetrization
operators \cite{DJ}.
Notice also that $A_3$ is a remarkable element of $H_n(q)$. Indeed,
$A_3=1+g_1+g_2+g_1g_2+g_2g_1+g_2g_1g_2$, and it is known that
the quotient of $H_n(q)$ by the ideal generated by $A_3$
is the Temperley-Lieb algebra $\text{TL}_n((2+q+q^{-1})^{-1})$ (see, e.g.,
\cite{GDJ}).

 We next show that for $q$ generic, the $A_k$'s are nonzero scalar
multiples
of idempotents.
We shall need the quantum factorials:
for $n\in{\Bbb N}$, 
$$1!_q:=1$$
$$n!_q:=(1+q+\dots+q^{n-1})(n-1)!_q.$$
\medskip

\noindent{\bf 2.1 Lemma} {\sl For $q\neq0, -1$ and for all $n\geq 2$ and
$j=1,\dots,n-1$, 
and $k=0,\dots,n-1$,
\begin{description}
\item{\rm a)} $g_jA_n=A_ng_j=qA_n,$
\item{\rm b)} $g_jA_{-n}=A_{-n}g_j=-A_{-n},$
\item{\rm c)} $\sigma^k(A_{n-k})A_{n}=(n-k)!_q A_{n},$
$\sigma^k(A_{-(n-k)})A_{-n}=(n-k)!_{1/q} A_{-n}.$
 In particular, $A_n^2=n!_qA_n$,
$A_{-n}^2=n!_{\frac{1}{q}}A_{-n}$.  
\item{\rm d)} If $q^n=1$, $A_n$ and $A_{-n}$ are central nilpotent
elements of $H_n(q)$, and therefore they belong to the radical of
that algebra.
\item{\rm e)} If $q^p\neq1$ for $p=2,\dots,n$, $E_n:=\frac{1}{n!_q}A_n$
and $E_{-n}:=\frac{1}{n!_{1/q}}A_{-n}$ are minimal and central idempotents
of $H_n(q)$.
\end{description}}\medskip

\noindent{\bf Proof} 
Notice that   $h_i:=-q^{-1}g_i$
satisfy the presentation relations for $H_\infty(q^{-1})$. 
Since the  $A_{-n}$  corresponds to
the  $A_n$ in the algebra $H_\infty(q^{-1})$ generated by the $h_i$'s,
 b) and  the second relation in c)  will follow from a) and the first
relation in c),
respectively. 
a)
We prove  that $g_jA_n=qA_n$ by induction on $n$. For $n=2$,
$A_2=1+g_1=(q+1)e_1$ and $g_1A_2=g_1+(g_1)^2=qg_1+q=q A_2.$
On the other hand for $i\geq1$,
$$g_1g_{i+1}g_i\dots g_1=g_{i+1}g_i\dots g_1g_2g_1=g_{i+1}\dots
g_2g_1g_2.$$
So, assuming that the claim holds  for $n$,
$$g_1g_{i+1}\dots g_1\sigma(A_n)=g_{i+1}\dots
g_1\sigma(g_1A_n)=qg_{i+1}\dots g_1\sigma(A_n),$$
thus
$$g_1A_{n+1}=g_1\sigma(A_n)+(g_1)^2\sigma(A_n)+q\sum_{i=2}^ng_i\dots
g_1\sigma(A_n)=qA_{n+1},$$
and the claim follows for $n+1$ and  $j=1$. 
Similar computations show that the claim holds also for 
 $j=2,\dots,n$. We now show that $A_ng_j=qA_n$ as well, for
$j=1,\dots,n-1$. First consider, for $n\in{\Bbb N}$,
the element $B_n$ defined by:
$B_1=1$,
$$B_{n+1}=\sigma(B_n)+\sigma(B_n)g_1+\dots+\sigma(B_n)g_1\dots g_n.$$
Similar arguments show that $B_ng_j=qB_n$ for $j=1,\dots,n-1$.
We claim that $B_n=A_n$ for all
$n$. 
It suffices to assume $n\geq1$.
In fact, $A_1=B_1=1$. Assume inductively that $A_k=B_k$ for $k=1,\dots,
n$. Then  computing $A_{n+1}$ by means of $A_n=B_n$ and $B_{n-1}$
gives
$$A_{n+1}=\sigma(A_n)+\sum_{i=1}^n g_i\dots
g_1\sigma(A_n)=$$
$$\sigma(A_n)+\sum_{i=1}^n g_i\dots
g_1\sigma(B_n)=$$
$$\sigma(A_n)+\sum_{i=1}^ng_i\dots
g_1\sigma^2(B_{n-1})+\sum_{i=1}^n\sum_{h=1}^{n-1}g_i\dots
g_1\sigma^2(B_{n-1})\sigma(g_1\dots g_h)=$$
$$\sigma(A_n)+\sum_{i=2}^n\sigma(g_{i-1}\dots
g_1\sigma(B_{n-1}))g_1+\sigma^2(B_{n-1})g_1+$$
$$\sum_{i=2}^n\sum_{h=1}^{n-1}\sigma(g_{i-1}\dots 
g_1\sigma(B_{n-1}))g_1\dots 
g_{h+1}+\sum_{h=1}^{n-1}\sigma^2(B_{n-1})g_1\dots 
g_{h+1}=$$
$$\sigma(B_n)+\sigma(B_n)g_1+\sum_{h=1}^{n-1}\sigma(B_n)g_1\dots
g_{h+1}=B_{n+1}.$$

c)  We fix $n$. $A_2=1+g_1$,  $\sigma^{n-1}(A_2)A_{n+1}=(1+q)A_{n+1}$ by a).
Taking into account the definition of the $A_j$'s, one obtains
iteratively
that $\sigma^k(A_{n+1-k})A_{n+1}=(n+1-k)!_q A_{n+1}$ for $k=n,\dots,0$.
d) By c), $A_n^2=A_{-n}^2=0$, so $A_n$ and $A_{-n}$ are nilpotent
in $H_\infty(q)$. On the other hand, they are central in $H_n(q)$,
therefore they are properly nilpotent elements of $H_n(q)$. e) It is
clear that $E_n$ and $E_{-n}$ are central idempotents of $H_n(q)$. a) and
b) show that  for all $X\in H_n(q)$, $E_nX$ is a scalar multiple of $E_n$,
so $E_n$ is a minimal idempotent.
\medskip

We are (eventually) interested in representing  Hecke algebras
in tensor $C^*$--categories. Therefore we look for
$^*$--involutions on $H_\infty(q)$.
If $q$ is real, a natural involution making $H_\infty(q)$  a
$^*$--algebra is that one for which
$g_i$'s become selfadjoint elements, 
 while, if $|q|=1$, 
one would like to  have 
$g_i^*=g_i^{-1}$. In both cases the $e_i$ become selfadjoint idempotents.

It is known that there is indeed a unique $^*$--involution on
$H_\infty(q)$ making 
the $e_i$'s  selfadjoint idempotents 
 if and only if
$q$ is  real or $|q|=1$. We shall refer to this involution as the
{\sl standard involution}.

\medskip

\noindent{\bf 2.2 Corollary} {\sl Assume that $q\neq0$, $q\neq-1$. For
$q$ 
 real 
or   $|q|=1$,
let us endow $H_\infty(q)$ with its standard involution. 
\begin{description}
\item{\rm a)}
If $q^n=1$, ${A_n}^*A_n=0$. In particular, $A_n$ lies in the kernel of any
Hilbert space $^*$--representation of $H_n(q)$.
\item{\rm b)}
If $q^p\neq1$
for $p=2,\dots,n$, the idempotents $E_n$ and $E_{-n}$ defined in the
previous 
lemma are selfadjoint.
\end{description}
}\medskip

\noindent{\bf Proof} Since the $e_i$ are selfadjoint, a straightfoward computation shows
that ${g_j}^*=\frac{\bar{q}+1}{q+1}g_j+\frac{\bar{q}-q}{q+1}$. Thus by 
a) in the previous lemma, ${g_j}^*A_n=\bar{q}A_n$ for $j=1,\dots,n-1$.
It follows that ${A_n}^*A_n=(1+\bar{q}+\dots+(\bar{q})^{n-1})
\sigma(A_{n-1}^*)A_n$. If $q^n=1$, ${A_n}^*A_n=0$, and this shows a). b) 
If instead $q^p\neq1$,
$p=2,\dots,n$,  ${A_n}^*A_n=n!_{\bar q}A_n=|n!_q|^2E_n$, which shows
that
$E_n$ is selfadjoint. Replacing the $g_j$'s by the  $-q^{-1}g_j$'s and 
$q$ by $q^{-1}$, amounts to showing that $E_{-n}$ is selfadjoint as
well.\medskip

For $q\neq0$, $q\neq-1$, $q^p\neq1$, $p=2,\dots, n$, $H_n(q)$ is a
semisimple
algebra \cite{We}: as for the the symmetric group on $n$ symbols, its
irreducibles $\{\pi_\lambda\}$ are labeled by Young diagrams 
$\lambda=[\lambda_1,\dots,\lambda_k]$ 
with $n$
squares. For the  explicit formula of the representation 
$\pi_\lambda$ we shall apply  formula  
 $(2.3)$ in \cite{We} to our idempotents $e_i$'s. Notice, however, that 
we adopt a different convention than in 
\cite{We}:  our idempotents $e_i$ correspond to Wenzl's $1-e_i$'s. 

Let us assume $q$ real, $q\neq-1$, or $|q|=1$, and let us endow
$H_\infty(q)$ with its standard $^*$--involution.
A $^*$--representation of $H_n(q)$, with $n$ possibly
infinite, on a Hilbert
space is called {\sl trivial} if it is the direct sum of
the one dimensional representations: $\pi_0(e_i)=0$, $\pi_1(e_i)=1$
for all $i$.  
$H_n(q)$  admits a non trivial
Hilbert space $^*$--representation for $n$ arbitrarily large 
 only if $q\geq 0$ or $q=e^{\frac{2\pi i}{m}}$, for some $m\in{\Bbb Z}$
with
$|m|=3,4,\dots$. 
For $q>0$ the above representations $\pi_\lambda$ are
$^*$--representations.

For   $q=e^{\frac{2\pi i}{m}}$,
with $m\in{\Bbb Z}$, $|m|\geq4$, and for any positive integer
$k\leq m-1$, Wenzl
defines
a  semisimple, irreducible
$^*$--representation $\pi_\lambda^{(k,m)}$ of $H_n(q)$
associated with every $(k,m)$-diagram 
$\lambda=[\lambda_1,\dots,\lambda_k]$ 
(a Young diagram with $n$ squares,
at most $k$ rows and such that 
$\lambda_1-\lambda_k\leq
m-k$).
Let $\Lambda_n^{(m)}$ denote the collection of all Young $(k,m)$-diagrams
with $n$ squares, for some $k\leq m-1$. The semisimple 
representation
$\pi_n^{(m)}=\oplus_{\lambda\in\Lambda_n^{(m)}}\pi_\lambda$ of $H_n(q)$
is, in general, not faithful.
\medskip

\noindent{\bf 2.3 Lemma} {\sl Assume that $q\neq0,-1$, $q^p\neq1$,
$p=2,\dots,k$.
For
any positive integer $k$, $E_k$ (resp. $E_{-k}$ ) is the minimal
central idempotent of $H_k(q)$ corresponding to the Young diagram $[1^k]$
(resp. $[k]$).}\medskip

\noindent{\bf Proof} By Lemma 2.1 e) $E_k$ is a minimal central  
idempotent 
of $H_k(q)$. It defines a one dimensional representation
$\chi$ of $H_k(q)$ by $XE_k=\chi(X)E_k$ and such that
$\chi(g_i)=q$, or, equivalently, $\chi(e_i)=1$, $i\leq k-1$.
We are left to show that $\pi_{[1^k]}(e_i)=1$, and this follows
from formula $(2.3)$ in \cite{We}. 
\medskip

For generic values of $q$, $H_n(q)$ is a semisimple algebra.
Its irreducibles are labeled  by Young diagrams with $n$ squares.
We have the problem of determining the Young diagrams  corresponding to
those irreducibles whose central supports sum up to the central
support of the ideal of $H_n(q)$
generated by $E_k$, for
$k\leq n$. Notice that,
for $k=3$, the quotient by that ideal is precisely the Temperley--Lieb
algebra
$\text{TL}_n((q+q^{-1}+2)^{-1})$,
which is known to correspond, in the above sense, to the set of diagrams
with at most $2$ rows,
see \cite{GHJ}.\medskip

\noindent{\bf 2.4 Proposition} {\sl If $q\neq0$, $q\neq-1$, $q^p\neq1$,
$p=2,\dots,n$, 
for any positive
integer $k\leq n$, the ideal of
$H_n(q)$ generated by the idempotent $E_k$ (resp. $E_{-k}$) corresponds to
the set of
Young diagrams with at least $k$ rows (resp. columns).}\medskip

\noindent{\bf Proof}
Let $p_\lambda$ denote a minimal central idempotent
of $H_n(q)$ corresponding to the Young diagram $\lambda$. By the
previous lemma $E_k$ is
  the minimal central idempotent of $H_k(q)$ corresponding 
to the Young diagram $[1^k]$.
In order 
to decide whether or not $p_\lambda E_k\neq0$ it is necessary and
sufficient  that $\lambda>[1^k]$ (see formula $(2.6)$ in \cite{We}) i.e.
that
$\lambda$ contain at least 
$k$ rows.
One similarly shows the remaining part.

\bigskip
\end{section}

\begin{section}{Hecke  symmetries in tensor $C^*$--categories}

Let us introduce {\sl the braid category} ${\Bbb B}$ with set of objects
the non
negative 
integers and arrows 
$$(n,m)=0,\quad \text{if } n\neq m,$$
$$(0,0)=(1,1)={\Bbb C},$$
$$(n,n)={\Bbb C}[{\Bbb B}_n],\quad n\geq 2,$$
the complex group algebra of the braid group ${\Bbb B}_n$ on $n-1$
generators: $g_1,\dots,g_{n-1}$.
with their natural structure of complex vector spaces, and composition
of arrows arising from the algebra structure of $(n,n)$.
${\Bbb B}$ becomes a strict tensor category if 
we define the tensor product between objects as $n\otimes m=n+m$
and between arrows as
$$S\otimes T=S\sigma^n(T),\quad S\in(n,n), T\in(m,m),$$

where $\sigma$ denotes, as before, the monomorphism ${\Bbb B}_m\to{\Bbb
B}_{m+1}$
taking $g_i$ to $g_{i+1}$ for $i=1,\dots, m-1$.

Let us consider a tensor category ${\cal T}$ with objects the tensor
powers of a 
single object $\rho$, with  tensor identity object given by $\iota$.

A {\sl braided  symmetry} for $\rho$ will be 
given by a tensor functor
$$\varepsilon: {\Bbb B}\to {\cal T}$$
such that $\varepsilon(0)=\iota$ and $\varepsilon(1)=\rho$.

 More explicitly, we shall need group
representations
$\varepsilon_n: {\Bbb B}_n\to(\rho^n,\rho^n)$
such that, for $b\in {\Bbb B}_n,$ 
$$\varepsilon_{n+1}(b)=\varepsilon_n(b)\otimes 1_\rho,$$
$$\varepsilon_{n+1}(\sigma(b))=1_\rho\otimes \varepsilon_n(b).$$

It is convenient to select 
the  subcategory ${\cal T}^\varepsilon$
of ${\cal T}$ with the same objects as ${\cal T}$ and arrows
$$(\rho^n,\rho^m)^\varepsilon=\{T\in(\rho^n,\rho^m):
\varepsilon(g_1\dots g_m)T\otimes
1_\rho=1_\rho\otimes T\varepsilon(g_1\dots g_n)\}.
\eqno(3.1)$$

${\cal T}$ is called a {\sl braided tensor category} if ${\cal
T}^\varepsilon={\cal T}$.

In a similar way, one can define {\sl the Hecke category} $H(q)$
with objects the non negative integers and arrows $(n,m)=0$ for $n\neq m$
and $(n,n)=H_n(q)$.

A braided symmetry $\varepsilon:{\Bbb B}\to{\cal T}$
will be called
{\sl a Hecke  symmetry} if 
it factors through a functor
$\varepsilon: H(q)\to{\cal T}$ from the Hecke category.

If ${\cal T}$ equals ${\Bbb B}$ (resp. $H(q)$), the
identity functor
does define a braided symmetry (resp. a Hecke  symmetry) for the
object $1$
making it into a braided tensor category.
The only nontrivial relation that 
we need to check is $(3.1)$, namely that
$$g_1\dots g_nb=\sigma(b)g_1\dots g_n,\quad b\in {\Bbb B}_n,$$
which follow from the presentation relations $(2.1)$--$(2.2)$ of the braid
group.

Let us assume that $H(q)$ admits a $^*$--involution making it into
a tensor $^*$--category.
In this case, a  $^*$--preserving Hecke  symmetry into another
tensor $^*$--category ${\cal T}$ will be called a {\sl 
$^*$--symmetry}.

As an example, if
 $q$ is real or if $|q|=1$, $H(q)$ becomes naturally a tensor
$^*$--category
with the involution inherited from the standard involution of 
$H_n(q)$.

Let $\Phi$ be a left inverse of $\rho$:
a set $\{\Phi_n, n\in{\Bbb N}\}$ of $(\iota,\iota)$--linear mappings 
$$\Phi_n:(\rho^n,\rho^n)\to (\rho^{n-1},\rho^{n-1})$$
(with  $\rho^0=\iota$) preserving
right tensoring by $1_\rho$ and satisfying for $S\in(\rho^n,\rho^n)$, 
$T\in(\rho^{n-1},\rho^{n-1}),$
$$\Phi_n(S1_\rho\otimes T)=\Phi_n(S)T,$$
$$\Phi_n(1_\rho\otimes TS)=T\Phi_n(S).$$
If ${\cal T}$ is a tensor $C^*$--category, $\Phi$ will be called
positive if each $\Phi_n$ is positive. The bimodule property then implies
that $\Phi_n$ is completely positive.

\medskip

\noindent{\bf 3.1 Lemma} {\sl Let $\varepsilon$ be a braided symmetry
for an object $\rho$ of a tensor $C^*$--category ${\cal T}$,
and let $\Phi$ be a  positive left inverse of $\rho$ such that
$\Phi_2(\varepsilon(g_1))$ is invertible. 
If $T\in(\rho^n,\rho^m)^\varepsilon$ satisfies $\Phi_n(T^*T)=0$
then $T=0$.}\medskip

\noindent{\bf Proof} Set
$\varepsilon_n:=\varepsilon(g_1\dots g_n)$. 
Since $\varepsilon_m T\otimes 1_\rho=1_\rho\otimes
T\varepsilon_n$,
 $T^*\otimes 1_\rho{\varepsilon_m}^*\varepsilon_m
T\otimes1_\rho=\varepsilon_n^*1_\rho\otimes T^*T\varepsilon_n$.
Assume that $\Phi_n(T^*T)=0$. Then  the left hand side of the
above
equation  is in the kernel of $\Phi_{n+1}$ 
by positivity. Writing
$\varepsilon_n=\varepsilon(g_1)\otimes
1_{\rho^{n-1}}1_\rho\otimes\varepsilon_{n-1}$ shows that
$\varepsilon(g_1)^*\otimes1_{\rho^{n-1}}1_\rho\otimes
T^*T\varepsilon(g_1)\otimes1_{\rho^{n-1}}$ lies in the kernel
of $\Phi_{n+1}$. Complete positivity implies 
$\Phi_{n+1}(A^*A)\geq\Phi_{n+1}(A)^*\Phi_{n+1}(A)$, which, in turn,
implies $0=\Phi_{n+1}(1_\rho\otimes
T\varepsilon(g_1)\otimes1_{\rho^{n-1}})=
T\Phi_2(\varepsilon(g_1))\otimes1_{\rho^{n-1}}$,
and the conclusion follows.
\medskip

For a given Hecke  $^*$--symmetry $\varepsilon: H(q)\to{\cal
T}_\rho$ we have the problem of determining the kernels of 
the corresponding $^*$--homomorphisms $H_n(q)\to(\rho^n,\rho^n)$.
If $\rho$ has a positive left inverse $\Phi$ such that
$\Phi(\varepsilon(g_1))$ is an invertible element in $(\iota,\iota)$, we
show that this element determines
the kernels of $\varepsilon$ uniquely. 
\medskip

\noindent{\bf 3.2 Corollary} {\sl Let $H(q)$ be endowed with a
$^*$--involution making it into a tensor $^*$--category, and  let
$\varepsilon: H(q)\to{\cal
T}_\rho$ be a Hecke  $^*$--symmetry into a tensor $C^*$--category.
If there is a
positive left
inverse $\Phi$ for $\rho$ such that 
$\Phi_2(\varepsilon(g_1))$ is an invertible element in $(\iota,\iota)$
then, for all $n\in{\Bbb N}$, the kernel of the $^*$--representation 
$\varepsilon: H_n(q)\to(\rho^n,\rho^n)$ depends only on
$\Phi_2(\varepsilon(g_1))$.
}\medskip

\noindent{\bf Proof}
Let us consider the 
quotient map 
 $\pi_q:{\Bbb C}[{\Bbb B}_n]\to H_n(q)$. Notice that
$\pi_1: {\Bbb C}[{\Bbb B}_n]\to {\Bbb C}[{\Bbb P}_n]$
associates to a braid on $n$ threads
a corresponding permutation of $(1,\dots, n)$.
 There is a natural section
$s_q$ of $\pi_q$
(i.e. a  linear map $s_q:H_n(q)\to{\Bbb C}[{\Bbb B}_n]$ 
such that $\pi_q\circ s_q$ is the 
identity map) constructed in the following way. First we consider the 
subset $B_n$ of ${\Bbb B}_n$ defined by:
$B_1=\{1\}$, 
$$B_{n+1}=\sigma(B_n)\cup g_1\sigma(B_n)\cup\dots\cup g_n\dots
g_1\sigma(B_n).$$
Set
 $H_n=\pi_q(B_n)$, for all $n$.
$H_n$ is a linear basis of $H_n(q)$,
therefore the 
inverse map $s_q$ of the restriction of $\pi_q$ to  $\pi_q: B_n\to H_n$
extends uniquely to the desired linear section $s_q$. 
If $h\in H_n\subset H_{n}(q)$ is of the form
$h=\sigma(h')$, with $h'\in H_{n-1}$ then
$\Phi(\varepsilon(h))=\varepsilon(h')$.
If instead $h=g_i\dots g_2g_1\sigma(h'')$ with $h''\in H_{n-1}$,
then $\Phi(\varepsilon(h))=\lambda\varepsilon(h')$, where 
$h'=g_{i-1}\dots g_1h''$. In both cases, $h'$ is the result of the image under
$\pi_q\circ s_1$ of that permutation of ${\Bbb P}_{n-1}$ obtained from 
$\pi_1s_q(h)$ by deleting $1$ from its cycle in its decomposition into 
disjoint cycles,
and then writing $n-1$ in place of $n$. For all $n$, if  $h\in H_n$, 
$\Phi_2\circ\dots\circ\Phi_n(\varepsilon(h))$ is a power of 
$\lambda:=\Phi(\varepsilon(g_1))$, therefore we get
a  $(\iota,\iota)$--valued, positive linear map $\omega_\lambda$ on
$H_\infty(q)$, which, by
the above, can be computed explicitly: $\omega_\lambda(h)$  is the product
of  factors of the form $\lambda^{k-1}$  for each cycle on length $k$
in the decomposition
of $\pi_1s_q(h)$. Now notice that
$\varepsilon(H_n(q))\subset(\rho^n,\rho^n)^\varepsilon$, so, if
 $\lambda$
is invertible, by the previous lemma,
an element $h\in H_n(q)$ is in the kernel of $\varepsilon$ if and 
only if 
$\Phi_n(\varepsilon(h^*h))=0$. 
On the other hand, 
$\Phi_n(\varepsilon(H_n(q)))$ is contained in the algebra generated by
 $\varepsilon(H_{n-1}(q))$ and $(\iota, \iota)$, which is in turn 
contained in $(\rho^{n-1},\rho^{n-1})^\varepsilon$, therefore
$\Phi_n(\varepsilon(h^*h))=0$ is equivalent to
$\Phi_{n-1}(\Phi_n(\varepsilon(h^*h)))=0$.
Iterating the left inverse $n$ times we get that
$\varepsilon(h)=0$ if and only if  
$\omega_\lambda(h^*h)=0$.
\medskip

Recall that
if $g_i$  generate $H_\infty(q)$, $h_i=-\frac{g_i}{q}$ generate
$H_\infty(\frac{1}{q})$. We thus have an invertible tensor $^*$--functor
$\phi: H(\frac{1}{q})\to H(q)$ such that
$\phi(h_i)=-\frac{g_i}{q}$. If $\varepsilon$ is a $^*$--symmetry of
 $H(q)$ in ${\cal T}$, for $q\geq1$ or $q=e^{2\pi i/m}$,
$m=3,4,\dots$,
then $\varepsilon\circ\phi$ is  a $^*$--symmetry of
$H(\frac{1}{q})$ in the same category,
and any such $^*$--symmetry arises in this way. 
Therefore we can
 restrict the 
values of 
$q$ to $[1,+\infty)\cup \{e^{\frac{2\pi i}{m}},m=3,\dots\}$.  
\medskip

\noindent{\bf 3.3 Theorem} {\sl Assume that
$q\in[1,+\infty)\cup\{e^{\frac{2\pi i}{m}}; m=4,\dots\}$. Let us
endow $H(q)$ with its standard  $^*$--involution: $e_i^*=e_i$ for all $i$.
Let
$\varepsilon:
(H(q), ^*)\to{\cal T}$ be a
 Hecke $^*$--symmetry
for an object $\rho$ of a tensor $C^*$--category, and let $\Phi$ be a unital,
positive left
inverse of $\rho$ 
for which $\Phi(\varepsilon(g_1))=:\lambda$ is a complex number. 
Set, for $d\in{\Bbb Z}$, $d\neq0$, $\lambda_d:=
\frac{q^d(q-1)}{q^d-1}
$.
Then $\lambda$ satisfies one of the following conditions:
 \begin{description}
\item{\rm a)} if $q\geq 1$,  then either $0\leq\lambda\leq q-1$
or $\lambda=\lambda_d$ for some $d\in{\Bbb Z}$, $d\neq 0$.
If $\lambda\in[0,q-1]$, all the maps
$\varepsilon:H_n(q)\to(\rho^n,\rho^n)$ are
 monomorphisms,
whereas if $\lambda=\lambda_d$, the kernel of the homomorphism
$H_n(q)\to(\rho^n,\rho^n)$ is the ideal generated by
$E_{-d-\frac{|d|}{d}}.$
\item{\rm b)} if $q=e^{\frac{2\pi i}{m}}$ then 
$\lambda=\lambda_d$ for some 
$d\in{\Bbb Z}$, $|d|\leq m-1$, $d\neq 0$. The kernel
of $\varepsilon: H_n(q)\to(\rho^n,\rho^n)$
coincides, for  $d$ negative, with the kernel of  Wenzl's
representation $\pi_n^{(-d,m)}$,
and, for $d$ positive with the kernel of $\pi_n^{(d,m)}\circ\alpha$,
with $\alpha$ the automorphism of $H_n(q)$ defined by
$\alpha(g_i)=q-1-g_i$.
\end{description}}\medskip

\noindent{\bf Proof} In the proof of   corollary 2.2 we have seen
that $A_{n+1}^*A_{n+1}=(n+1)!_{\bar{q}}A_{n+1}$ and 
$A_{-(n+1)}^*A_{-(n+1)}=(n+1)!_qA_{-(n+1)}$, so
 $$\Phi(\varepsilon(A_{n+1}^*A_{n+1}))=
(n+1)!_{\bar{q}}\Phi(\varepsilon(A_{n+1}))=$$
$$(n+1)!_{\bar{q}}\varepsilon(A_n+\lambda A_n+\lambda g_1A_n+
\dots+\lambda g_{n-1}\dots g_1A_n)=$$
$$(n+1)!_{\bar{q}}[1+\lambda(1+ q+\dots+ q^{n-1})]\varepsilon(A_n)=$$
$$[1+\bar{q}+\dots+\bar{q}^n][1+\lambda(1+ q+\dots+
q^{n-1})]\varepsilon(A_n^*A_n).$$
Similarly,
$$\Phi(\varepsilon(A_{-(n+1)}^*A_{-(n+1)}))=$$
$$[1+(\bar{q})^{-1}+\dots+(\bar{q})^{-n}][1-q^{-1}\lambda(1+q^{-1}+\dots+
(q^{-1})^{n-1})]
\varepsilon(A_{-n}^*A_{-n}).$$
Assume $q\geq1$. If 
$\Phi(\varepsilon(A_k^*A_k))\neq0$  for all integers $k$,
 positivity of $\Phi$
implies $\lambda\in[0,q-1]$. If, on the contrary, 
$d$ is the integer with smallest
absolute value for which $\Phi(\varepsilon(A_d^*A_d)=0$, then necessarily 
$|d|\geq2$. If $d$ is negative,
$\lambda=\lambda_{-d-1}$, whereas if $d$ is positive,
$\lambda=\lambda_{-d+1}$.
In the case where $q=e^{2\pi i/m}$, $A_m^*A_m=0=A_{-m}^*A_{-m}$. 
Defining $d$ still as above, now implies
 $|d|\leq m$ and $d$ still needs to assume the stated values.

Let $\omega_\lambda$ be defined as in the proof of the previous corollary.
For $q\geq1$ and $\lambda\in[0,q-1]$, the previous computations show that
$\omega_\lambda({A_k}^*A_k)>0$ for all $k$, so
$\varepsilon({A_k}^*A_k)\neq0$, and this shows that $\varepsilon$ is
faithful. It is also evident that if $\lambda=\lambda_d$ then
$E_{-d-\frac{|d|}{d}}$, and therefore the ideal it generates in $H_n(q)$, 
lies in the kernel of $\varepsilon$. Taking into account the previous
corollary, the proof of a) will be complete if we  produce
 for each $d\in{\Bbb Z}$,
$d\neq 0$,
a  model   $^*$--symmetry $\varepsilon$
of $H_\infty(q)$,
with left inverse $\Phi$ for which $\Phi(\varepsilon(g_1))=\lambda_d$
and having the kernel as stated.
This is the content of the next section.
In case b),
let us assume that
$\Phi(\varepsilon(g_1))=\lambda_{-d}$,
for an appropriate positive $d$.
We have seen that an element $h\in H_\infty(q)$ is annihilated by 
$\varepsilon$ if and only if $\omega_{\lambda_{-d}}(h^*h)=0$.
Now $\omega_{\lambda_d}$ is a positive Markov trace on $H_\infty(q)$
factoring through the $C^*$--representation $\varepsilon$. This condition
says that $h$ lies in the kernel of the GNS representation associated to
this trace. By 
Theorem 3.6 in
\cite{We}, 
the latter is equivalent to
$\pi_n^{(d,m)}$.
\medskip

\end{section}

\begin{section}{The Model Hecke  $^*$--symmetry}

Notice that if $g_i$, $i=1,2,\dots,$ are generators 
of $H_\infty(q)$, then the elements $l_i:=q-1-g_i$, $i=1,2,\dots$
still satisfy the presentation relations of $H_\infty(q)$. In other words,
we have an  invertible tensor endofunctor $\alpha$ of
$H(q)$ 
which associates $l_i$ to $g_i$.

With every  Hecke symmetry $g_i\to\varepsilon(g_i)$ of
$H(q)$ in a tensor category we can associate 
 another Hecke symmetry of 
$H(q)$ in the same category by
$$\varepsilon'=\varepsilon\circ\alpha.$$ We shall refer to
$\varepsilon'$
as the {\sl dual} symmetry. 

If 
we endow
$H(q)$ with its standard
involution for appropriate values of $q$, $\alpha$ becomes a
$^*$-functor. 
Obviously, if  $\varepsilon$ is a $^*$--symmetry with respect to the
standard involution, 
$\varepsilon'$ is a $^*$--symmetry, as well.


Let $H$ be a finite dimensional vector space of dimension $d$. We consider 
the representation of the Hecke algebra $H_n(q)$ on $H^{\otimes n}$
discovered by Jimbo \cite{Ji} and Woronowicz \cite{Wo}. 
Let $g_q$ denote the operator on 
$H\otimes H$:
$$g_q\psi_i\otimes
\psi_j=-\mu\psi_j\otimes\psi_i,\quad i<j,$$
$$g_q\psi_i\otimes \psi_i=-\psi_i\otimes\psi_i$$
$$g_q\psi_i\otimes\psi_j=(q-1)\psi_i\otimes\psi_j-\mu\psi_j\otimes\psi_i,
\quad i>j,$$
with $\psi_1,\dots,\psi_d$ a basis of $H$ and $\mu$ a fixed
complex square root of $q$.
One can check that, on the vector space $H^{\otimes n}$, the operators 
$g_1=g_q\otimes 1_{H^{\otimes n-2}}$, $g_2=1_H\otimes g_q\otimes
1_{H^{\otimes n-3}},\dots$, $g_{n-1}=1_{H^{\otimes n-2}}\otimes g_q$
do satisfy all the relations $(2.1)$--$(2.3)$, and therefore
we obtain in this way a  model Hecke symmetry
$$\varepsilon: H(q)\to{\cal L}_H$$ with ${\cal L}_H$ the tensor
category of vector spaces with objects the tensor powers of $H$.

One can compute
the images of the   idempotents $e_i=\frac{1+g_i}{q+1}$
by  taking the translates
of 
$e_q:=\frac{1+g_q}{q+1}$,
$$e_q\psi_i\otimes \psi_j=\frac{1}{q+1}
(\psi_i\otimes\psi_j-\mu\psi_j\otimes\psi_i), \quad i<j$$
$$e_q\psi_i\otimes \psi_i=0$$
$$e_q\psi_i\otimes\psi_j=\frac{1}{q+1}
(q\psi_i\otimes\psi_j-\mu\psi_j\otimes\psi_i),
\quad i>j.$$

Let us consider first the special case where $d=2$.
It is then  known, and easy to check,  that the corresponding
representation of $H_n(q)$ is not
faithful for $n>2$:
it factors through a representation of the Temperley--Lieb algebra
$\text{TL}_n((2+q+q^{-1})^{-1})$, which is precisely the quotient
of the Hecke algebra $H_n(q)$ by the ideal generated by the idempotent
$E_3$. This fact generalizes to higher dimensions.\medskip

In general, we want to identify the kernels of
the   homomorphisms
$$\varepsilon: H_n(q)\to (H^{\otimes n}, H^{\otimes n}).$$
For values of $q$ for which the algebras $H_n(q)$ are semisimple, 
those kernels are
determined by specifying 
the quasiequivalence class of the above representation of $H_n(q)$
on  $H^{\otimes n}$.
\medskip

\noindent{\bf 4.1 Proposition} {\sl If $q\neq0$, $q\neq-1$, and $q^p\neq1$
for
$p=2,\dots,n$, the quasiequivalence class of $\varepsilon: H_n(q)\to
B(H^{\otimes n})$,
with $d\leq
n-1$,
 corresponds 
to the set of Young diagrams with $n$ squares and at most $d$ rows.
Alternatively, it coincides with the ideal of $H_n(q)$ generated by 
$E_{d+1}$. If $d\geq n$, $\varepsilon$ is faithful.}\medskip

\noindent{\bf Proof} The first assertion is shown in subsection 2.1 of
\cite{GW}. Notice that our model representation of $H_n(q)$ 
defined by $g_q$ is equivalent to the representation defined by
the operator $1-q-R$, where $R$ is defined in section 2 of \cite{GW}. 
Also, our irreducible representation $\pi_\lambda$ associated to the
Young diagram $\lambda$ corresponds to their $V_{\lambda'}$, where
$\lambda'$ is the diagram obtained from $\lambda$ exchanging rows with
columns.

In particular, if $d\geq n$, any Young diagram
with $n$ boxes gives rise to a subrepresentation of $\varepsilon$, so
$\varepsilon$ is faithful.

Therefore the kernel of $H_n(q)\to (H^{\otimes n}, H^{\otimes n})$
corresponds to the Young diagrams with $\geq d+1$ rows. The conclusion
follows from Prop 2.4.
\medskip

It follows that the dimension of $H$ is uniquely determined by the 
quasiequivalence class of the model symmetry $\varepsilon$.

An easy inductive argument shows that, for $k\in{\Bbb
N}$, $k\neq 0$,
$$\alpha(A_k)=q^{k(k-1)/2}A_{-k},$$
(see, e.g., \cite{Murphy95}).
Therefore,
if we have a  Hecke symmetry of $H(q)$ in a tensor category 
 with kernel the ideal generated by 
$E_{d+1}$, passing to the dual symmetry will give a Hecke symmetry 
with kernel the ideal generated by $E_{-d-1}$.
\medskip

In analogy with the case of a permutation symmetry treated in
\cite{DRduals}, we
introduce the concept of the integral dimension of a Hecke symmetry.
\medskip

\noindent{\bf 4.2 Definition} {\sl Assume that $q\neq0$, $q\neq-1$,
$q^n\neq1$
for all $n\geq2$.
Let $d$ be any integer with $|d|\geq2$.
We shall say that a Hecke symmetry $\varepsilon: H(q)\to{\cal
T}_\rho$ {\sl has dimension $d$} if, for all $n$,
the kernel of the homomorphism $H_n(q)\to(\rho^n,\rho^n)$
is the ideal generated by $E_{d+\frac{|d|}{d}}$.}\medskip 

In particular, for generic values of $q$, a $2$-dimensional 
Hecke symmetry is just a symmetry factoring over
the Temperley--Lieb algebras $\text{TL}_n((2+q+q^{-1})^{-1})$
and faithful on them.

Assume that $H$ is a Hilbert space and that $\psi_1,\dots,\psi_d$ is an
orthonormal basis.
The Hilbert space adjoint of $e_q$ is given by
$${e_q}^*\psi_i\otimes \psi_j=\frac{1}{\bar{q}+1}
(\psi_i\otimes\psi_j-\bar{\mu}\psi_j\otimes\psi_i), \quad i<j$$
$${e_q}^*\psi_i\otimes \psi_i=0$$
$${e_q}^*\psi_i\otimes\psi_j=\frac{1}{\bar{q}+1}
(\bar{q}\psi_i\otimes\psi_j-\bar{\mu}\psi_j\otimes\psi_i),
\quad i>j.$$
These computations show the following fact.

\medskip

\noindent{\bf 4.3 Proposition} {\sl
The idempotent  $e_q$ is selfadjoint if and only if $q$ is a positive real
number. Therefore
the Hecke symmetry
$\varepsilon$ is a $^*$--symmetry with respect to the standard
involution of $H(q)$ if and only if $q>0$.}\medskip

We next show by simple computations that 
in the case where $q$ is not a positive real, the Hilbert space $H$ can
not accomodate interesting $^*$--symmetries.
We introduce the symmetry
$\varepsilon_*:H(\overline{q})\to{\cal L}_H,$
which takes the basic generator of $H(\overline{q})$ 
to $g_q^*$, the Hilbert space adjoint of $g_q$.
\medskip

The following simple proposition shows that
if $|q|=1$,
the smallest
tensor
$^*$--subcategory of ${\cal L}_H$ containing the 
operators $g_q$ {\sl is permutation symmetric!}.\medskip

\noindent{\bf 4.4 Proposition} {\sl If $|q|=1$, the
smallest
$^*$--subalgebra
of $(H\otimes H, H\otimes H)$ containing $1$ and $g_q$ also contains 
 the permutation operator $\theta$ which exchanges the order
of factors in $H\otimes H$.}\medskip

\noindent{\bf Proof} We can assume $q\neq 1$. An  computation shows
that if $|q|=1$,
$$({g_q}^{-1})^*\psi_i\otimes
\psi_j=-\mu\psi_j\otimes\psi_i+(q-1)\psi_i\otimes\psi_j\quad i<j,$$
$$({g_q}^{-1})^*\psi_i\otimes \psi_i=-\psi_i\otimes\psi_i$$
$$({g_q}^{-1})^*\psi_i\otimes\psi_j=-\mu\psi_j\otimes\psi_i,
\quad i>j$$
therefore $({g_q}^{-1})^*-g_q=(q-1)P$,
where 
$$P\psi_i\otimes\psi_j=\psi_i\otimes\psi_j,\quad i<j,$$
$$P\psi_i\otimes\psi_i=0,$$
$$P\psi_i\otimes\psi_j=-\psi_i\otimes\psi_j,\quad i>j.$$
On the other hand,
$$1/2[({g_q}^{-1})^*g_q+g_q({g_q}^{-1})^*]-q=(1-q)T$$
with
$$T\psi_i\otimes\psi_j=\mu\psi_j\otimes\psi_i,\quad i\neq j,$$
$$T\psi_i\otimes\psi_i=\psi_i\otimes\psi_i.$$
Now the equality  $T=\mu\theta+(1-\mu)(1-P^2)$
completes the proof.
\bigskip

\end{section}


\begin{section} {Special Intertwiners and Conjugates}

In this section we specialize to the case where $q>0$.
Let $\mu$ denote the positive square
root of $q$. 
Consider the element of $H^{\otimes d}$ defined by
$$S=\sum_{p\in{\Bbb
P}_d}(-\mu)^{i(p)}\psi_{p(1)}\otimes\dots\otimes\psi_{p(d)}.$$

This is the fundamental intertwiner of the
representation
category of the Woronowicz quantum group $S_{\mu}U(d)$.
Here $i(p)$ is number of inversed pairs on $p$, i.e. 
the cardinality of the set $\{(i,j): i<j,  p(i)>p(j)\}$.
Note that $S$ is not normalized for $d>1$. In fact, 
$$\|S\|^2=\sum_{p\in{\Bbb P}_d}q^{i(p)}$$
and this is known to coincide with  the quantum $d$ factorial defined at
the
beginning of section $2$. 
Indeed, the claim is true for $d=1$. For $d>1$, let us represent ${\Bbb
P}_d$ as the disjoint union of left cosets 
${\Bbb P}_d={\Bbb P}'_{d-1}\cup\cup_{h=1}^{d-1}(h h+1)\dots (12) {\Bbb
P}'_{d-1}$, with ${\Bbb P}'_{d-1}:=\{p\in{\Bbb P}_d: p(1)=1\}$, which is a
copy of ${\Bbb P}_{d-1}$ contained in ${\Bbb P}_d$. Since $i((h
h+1)\dots(12)q)=h+i(q)$ for all
$q\in{\Bbb P}'_{d-1}$, 
$$\sum_{p\in{\Bbb P}_d}q^{i(p)}=(1+q+\dots+q^{d-1})\sum_{p\in{\Bbb
P}_{d-1}} q^{i(p)}.$$
More generally, choose
 $n\leq d$ increasing indices $i_1<i_2<\dots<i_n$ among $\{1,\dots,d\}$
and consider the $n$--dimensional 
subspace $H_{i_1,\dots, i_n}$ generated by $\psi_{i_1},\dots,\psi_{i_n}$
and the corresponding  symmetric tensor:
$$S_{i_1,\dots,i_n}:=\sum_{p\in{\Bbb
P}_n}(-\mu)^{i(p)}\psi_{i_{p(1)}}\otimes\dots\otimes\psi_{i_{p(n)}}.$$

\bigskip

\noindent{\bf 5.1 Lemma} {\sl Let $\varepsilon$ be the model Hecke
$^*$--symmetry defined in Sect. 4, and
$\varepsilon'=\varepsilon\circ\alpha$ its dual
symmetry.
 \begin{description}
\item{\rm a)}For $i=1,\dots,n-1,$
$$\varepsilon(g_i)S_{i_1,\dots,i_n}=qS_{i_1,\dots,i_n},$$
$$\varepsilon'(g_i)S_{i_1,\dots,i_n}=-S_{i_1,\dots,i_n},$$

\item{\rm b)} $\varepsilon(A_n)$ 
has support on the linear span of
vectors $\psi_{i_1}\otimes\dots\otimes\psi_{i_n}$ with $i_h\neq i_k$ for
$h\neq k$. In particular, $\varepsilon(A_n)=0$ for $n>d$. Furthermore,
for $n\leq d$, and for
 $i_1<i_2<\dots<i_n$,
$$\varepsilon(A_n)\psi_{i_{p^{-1}(1)}}\otimes\dots\otimes\psi_{i_{p^{-1}(n)}}=
(-\mu)^{i(p)}S_{i_1,\dots,i_n}.$$



\item{\rm c)} 
$$\sum_{i_1<\dots<i_n}S_{i_1,\dots,i_n}{{S}_{i_1,\dots,i_n}}^*=
\varepsilon(A_n).$$
In particular,
$$S{S}^*=\varepsilon(A_d).$$

\item{\rm d)} For $n\leq d$,  $i_1<i_2<\dots< i_n$, 
 $$\varepsilon(g_1\dots g_n)S_{i_1,\dots,i_n}\otimes 1_{H_{i_1,\dots, i_n}}=
-(-\mu)^{n-1}1_{H_{i_1,\dots, i_n}}\otimes S_{i_1,\dots,i_n},$$
 $$\varepsilon'(g_1\dots g_n)S_{i_1,\dots,i_n}\otimes
1_{H_{i_1,\dots, i_n}}=\mu^{n+1}1_{H_{i_1,\dots, i_n}}\otimes S_{i_1,\dots,i_n}.$$
In particular,
 $$\varepsilon(g_1\dots g_d)(S\otimes
1_H)=-(-\mu)^{d-1}1_H\otimes S,$$
 $$\varepsilon'(g_1\dots g_d)(S\otimes
1_H)=\mu^{d+1}1_H\otimes S.$$

\item{\rm e)} Set: $H^{(0)}=\text{span}\{\psi_r, r>i_n\}$, 
$H^{(n)}=\text{span}\{\psi_r, r<i_1\}$
if those sets are not empty, and,
for $h=1,\dots, n-1$,  
 $H^{(h)}=\text{span}\{\psi_r, i_{n-h}<r<i_{n-h+1}\}$, if $|i_{n-h+1}-i_{n-h}|\geq 2$.
Then
$$S_{i_1,\dots,i_n}^*\otimes 1_{H}
1_{H^{\otimes n-1}}\otimes
g_qS_{i_1,\dots,i_n}\otimes
1_{H}=-(n-1)!_q1_{H_{i_1,\dots,i_n}}+$$
$$(n-1)!_q\sum_h
(q^h-1)1_{H^{(h)}}$$
where the sum is taken over all $h=0,\dots n$ for which $H^{h}$ makes
sense.
\end{description}}
\medskip

\noindent{\bf Proof} 
We omit the proof of a), as it is straightforward.
b) Consider, for $n\in{\Bbb N}$ a vector in $H^{\otimes n}$
of the form $\psi=\psi_{i_1}\otimes\dots\otimes\psi_{i_n}$.
One can easily check by induction on $n$ that, if $i_1<i_2<\dots<i_n$,
$$\varepsilon(s(p))\psi=(-\mu)^{i(p)}
\psi_{i_{p^{-1}(1)}}\otimes\dots\otimes\psi_{i_{p^{-1}(n)}},$$ with 
$s:{\Bbb P}_n\to H_n(q)$ the section already considered in the proof of
Corollary 3.2.
Combining this  observation with the equation
$\varepsilon(A_n)\varepsilon(s(p))=q^{i(p)}\varepsilon(A_n)$, which
follows from a) in Lemma 2.1, shows that
$\varepsilon(A_n)\psi_{i_{p^{-1}(1)}}\otimes\dots\otimes\psi_{i_{p^{-1}(n)}}=
(-\mu)^{i(p)}\varepsilon(A_n)\psi$. Also,
$\varepsilon(A_n)\psi=\sum_{p\in{\Bbb P}_n}\varepsilon(s(p))\psi=
\sum_{p\in{\Bbb
P}_n}(-\mu)^{i(p)}\psi_{i_{p^{-1}(1)}}\otimes\dots\otimes\psi_{i_{p^{-1}(n)}}=
S_{i_1,\dots, i_n}$ as $i(p)=i(p^{-1})$.
If, more generally, $i_1\leq i_2\leq\dots\leq i_n$ and, for instance,
$i_1=i_2$ then applying both sides
of $\varepsilon(A_ng_1)=q\varepsilon(A_n)$ to $\psi$ yields
$\varepsilon(A_n)\psi=0$, as $q\neq-1$. On the other hand,
as above, the relation
$\varepsilon(A_n)\varepsilon(s(p))=q^{i(p)}\varepsilon(A_n)$, due to
Lemma 2.1 a), with 
$s$ the section defined as in the proof of Corollary 3.2, applied to
$\psi$ shows 
that $\varepsilon(A_n)
\psi_{i_{p^{-1}(1)}}\otimes\dots\otimes\psi_{i_{p^{-1}(n)}}$ is a scalar
multiple of $\varepsilon(A_n)\psi$, and the latter vanishes if at least
two 
indices repeat. 
c) Notice that if two ordered $n$-tuples $i_1<i_2<\dots<i_n$,
$j_1<\dots<j_n$ are different, the corresponding $S$'s are orthogonal.
Also,
$\sum_{j_1<\dots<j_n}S_{j_1,\dots,j_n}{{S}_{j_1,\dots,j_n}}^*$
has the same support as $\varepsilon(A_n)$
and, for $i_1<\dots<i_n$, 
$$\sum_{j_1<\dots<j_n}S_{j_1,\dots,j_n}{{S}_{j_1,\dots,j_n}}^*
\psi_{i_{p(1)}}\otimes\dots\otimes\psi_{i_{p(n)}}=
(-\mu)^{i(p)}S_{i_1,\dots,i_n}=$$
$$\varepsilon(A_n)\psi_{i_{p(1)}}\otimes\dots\otimes\psi_{i_{p(n)}}.$$
d)
 Since
$S_{i_1,\dots,i_n}=\varepsilon(A_n)\psi_{i_1}\otimes\dots\otimes\psi_{i_n}$, 
$$\varepsilon(g_1\dots g_n)S_{i_1,\dots,i_n}\otimes\psi_{i_n}=
\varepsilon(g_1\dots g_n
A_n)\psi_{i_1}\otimes\dots\otimes\psi_{i_n}\otimes\psi_{i_n}=$$
$$\varepsilon(\sigma(A_n)g_1\dots
g_n)\psi_{i_1}\otimes\dots\otimes\psi_{i_n}\otimes\psi_{i_n}=
-\varepsilon(\sigma(A_n))\varepsilon(g_1\dots g_{n-1})\psi_{i_1}\otimes
\dots\otimes\psi_{i_n}\otimes\psi_{i_n}=$$
$$-(-\mu)^{n-1}\varepsilon(\sigma(A_n))
\psi_{i_n}\otimes\psi_{i_1}\otimes\dots\otimes\psi_{i_n}=
-(-\mu)^{n-1}\psi_{i_n}\otimes S_{i_1,\dots,i_n},$$
while, for $i_n>r$ and $r\in\{i_1,\dots, i_n\}$,
$$\varepsilon(g_1\dots g_n)S_{i_1,\dots, i_n}\otimes\psi_r=
-\mu\varepsilon(\sigma(A_n))\varepsilon(g_1\dots
g_{n-1})\psi_{i_1}\otimes\dots\otimes\psi_{i_{n-1}}\otimes\psi_r\otimes\psi_{i_n}+$$
$$(q-1)\varepsilon(\sigma(A_n))\varepsilon(g_1\dots 
g_{n-1})\psi_{i_1}\otimes\dots\otimes\psi_{i_n}\otimes\psi_r.$$
Now $$\varepsilon(g_1\dots g_{n-1})\psi_{i_1}\otimes\dots\otimes
\psi_{i_n}\otimes\psi_r=
(-\mu)^{n-1}\psi_{i_n}\otimes\psi_{i_1}\otimes\dots\otimes\psi_{i_{n-1}}\otimes\psi_r$$ 
so the last addendum in the above sum vanishes as the index $r$ appears
twice in $\psi_{i_1}\otimes\dots\otimes\psi_{i_{n-1}}\otimes\psi_r$.
Iterating
this
procedure gives the desired result.

For the analogous relation relative to the dual symmetry, similar
computations will take to the desired result. One can proceed in the
following way. Let $p_n$ be that permutation that reverses the order
of the ordered integers $i_1<\dots<i_n$, and write
$$S_{i_1,\dots,i_n}=(-\mu)^{-i(p_n)}\varepsilon(A_n)\psi_{i_n}\otimes
\dots\psi_{i_1}.$$ Then 
$$\varepsilon'(g_1\dots
g_n)S_{i_1,\dots,i_n}\otimes\psi_r=
(-\mu)^{-i(p_n)}\varepsilon(\sigma(A_n))\varepsilon'(g_1\dots 
g_n)\psi_{i_n}\otimes\dots\psi_{i_1}\otimes\psi_r.$$
Now for $r=i_1$, as before, we easily get the result, and for $r>i_1$
we just need to know that for $i\neq j$, 
the new coefficient
$(\psi_i\otimes\psi_j,
\varepsilon'(g_1)\psi_j\otimes\psi_i)$ is  $\mu$.

We give another  proof  in the case $n=d$. By
c),
$$d!_q\varepsilon'(g_1\dots g_d)S\otimes 1_H=\varepsilon'(g_1\dots
g_d)\varepsilon(A_d)S\otimes 1_H=$$
$$\varepsilon(\sigma(A_d))\varepsilon'(g_1\dots
g_d)S\otimes 1_H=
1_H\otimes S(1_H\otimes S^*\varepsilon'(g_1\dots g_d)S\otimes 1_H).$$
Now by the previous part
$$1_H\otimes S^*\varepsilon'(g_1\dots
g_d)S\otimes 1_H=$$
$$(-1)^d\frac{1}{\mu^{d-1}}
S^*\otimes 1_H\varepsilon(g_d\dots g_1)\varepsilon'(g_1\dots
g_d)S\otimes 1_H=\mu^{d+1} d!_q
$$
as $\varepsilon(g_i)\varepsilon'(g_i)=-q$.


e) We apply  the left hand side  to all basis vectors $\psi_r$. We
start from the case where $r$ is of the form $i_j$. 
$$S_{i_1,\dots,i_n}^*\otimes1_{H}(
-\sum_{p(n)=j}
(-\mu)^{i(p)}
\psi_{i_{p(1)}}\otimes\dots\otimes\psi_{i_{p(n)}}\otimes
\psi_{i_j}+$$
$$-\mu\sum_{p(n)\neq j}(-\mu)^{i(p)}
\psi_{i_{p(1)}}\otimes\dots\otimes\psi_{i_{p(n-1)}}\otimes
\psi_{i_j}
\otimes\psi_{i_{p(n)}}+$$
$$(q-1)\sum_{p(n)>j}
(-\mu)^{i(p)}
\psi_{i_{p(1)}}\otimes\dots\otimes\psi_{i_{p(n)}}\otimes
\psi_{i_j}).$$
For any permutation $p$ for which $p(n)\neq j$, 
 $j$ appears twice in $(p(1),p(2),\dots,p(n-1), j)$, therefore
the second sum, when multiplied on the left by 
$S_{i_1,\dots,i_n}^*\otimes 1_{H_{i_1,\dots,i_n}}$, 
vanishes, and the computation equals
$$-\sum_{p(n)=j}
(-\mu)^{i(p)}
S_{i_1,\dots,i_n}^*\otimes1_{H}\psi_{i_{p(1)}}
\otimes\dots\otimes\psi_{i_{p(n)}}\otimes\psi_{i_j}+$$
$$(q-1)\sum_{k=j+1}^n
\sum_{p(n)=k}
(-\mu)^{i(p)}
S_{i_1,\dots,i_n}^*\otimes1_{H}
\psi_{i_{p(1)}}\otimes\dots\otimes\psi_{i_{p(n)}}\otimes
\psi_{i_j}.$$
Let us write each permutation for which $p(n)=k$ in the 
form $p=(k k-1)(k-1 k-2)\dots(2 1)p''(1 2 3\dots n-1 n)$ 
with $p''(1)=1$.
We claim that  $i(p)=i(p'')+n-k$. 
In fact, for any permutation $q$
for which $q(1)=h$ the permutation $(h-1 h)q$ has one 
less inversed 
pair than $q$, so $i((h-1 h)q)=i(q)-1$, which shows that 
$i(p'')=i(p(1\dots n)^{-1})-(k-1)$. 
Now if we want to compare the inversed
pairs of $p$ and $p(1\dots n)^{-1}$
we see that the two sets of inversed pairs 
have a common intersection with
cardinality, say $m$,
but  for $p$
we need to count how many inversed pairs we have in the set
$\{(p(1),k),\dots,(p(n-1),k)\}$, which are $n-k$, 
whereas for $p(1\dots
n)^{-1}$ we need the cardinality of the subset of inversed pairs in
$\{(k,p(1)),\dots,(k,p(n-1))\}$, which is $k-1$.
 Therefore 
$i(p)=m+n-k$ and $i(p(1\dots n)^{-1}=m+k-1$,
 which implies $i(p(1\dots
n)^{-1})=i(p)+2k-n-1$ and the claim is proved.
Thus
 $$\sum_{p(n)=k}
(-\mu)^{i(p)}
S_{i_1,\dots,i_n}^*\otimes1_{H}
\psi_{i_{p(1)}}\otimes\dots\otimes\psi_{i_{p(n)}}\otimes
\psi_{i_j}=
\sum_{p''\in{\Bbb P}_{n-1}}q^{i(p'')+n-k}\psi_{i_j},$$
and the conclusion follows from simple computations in this case.
Assume now that $r$ is none of the $\psi_{i_j}$'s, and, to fix ideas, assume that
$i_{n-h}<r<i_{n-h+1}$. Then similar computations yield:
$$S_{i_1,\dots, i_n}^*\otimes 1_H 1_{H^{\otimes(n-1)}}\otimes g_q 
S_{i_1,\dots,i_n}\otimes 1_H(\psi_r)=$$
$$S_{i_1,\dots, i_n}^*\otimes 1_H(-\mu\sum_{p}
(-\mu)^{i(p)}\psi_{i_{p(1)}}\otimes\dots\otimes\psi_{i_{p(n-1)}}\otimes\psi_r\otimes
\psi_{i_{p(n)}}+$$
$$(q-1)\sum_{i_{p(n)}>r}(-\mu)^{i(p)}\psi_{i_{p(1)}}\otimes\dots\otimes
\psi_{i_{p(n)}}\psi_r)=(q-1)\sum_{i_{p(n)}>r}q^{i(p)}\psi_r.$$
Notice that the condition $i_{p(n)}>r$ amounts to
$p(n)=k$ with $k=n-h+1,\dots,n$. Therefore, as before, 
 we write each such permutation in the form $p=(k k-1)\dots(21)p''(12\dots n)$ with
$p''(1)=1$ and $i(p)=i(p'')+n-k$. We thus conclude that the above term equals
$$(q-1)\sum_{k=n-h+1}^n\sum_{p''}q^{i(p'')+n-k}\psi_r=
(q-1)(n-1)!_q\sum_{k=n-h+1}^nq^k\psi_r=$$
$$(q^h-1)(n-1)!_q\psi_r.$$
In the case where $r<i_1$ the computations are simpler, as $i_{p(n)}>r$ for all $p\in{\Bbb P}_n$,
so the result is $(q-1)n!_q\psi_r=(q^n-1)(n-1)!_q\psi_r$. Finally, if $r>i_n$
then $S_{i_1,\dots, i_n}^*\otimes 1_H1_{H^{\otimes (n-1)}}
\otimes g_qS_{i_1,\dots, i_n}\otimes 1_H\psi_r=0.$

\medskip

The previous lemma has various consequences. We start 
with the following one.
\medskip

\noindent{\bf 5.2 Corollary} {\sl The representation category of $S_\mu
U(d)$
is a braided tensor category.}\medskip

\noindent{\bf Proof} Let $\omega$ be a complex number such that
$\omega^d=-(-\mu)^{d-1}$, and consider the representations
$\tilde{\varepsilon}_n:{\Bbb B }_n\to (H^{\otimes n}, H^{\otimes n})$
given by $\tilde{\varepsilon}_n(g_i)=\frac{\varepsilon(g_i)}{\omega}$.
This is obviously a braided symmetry in the category ${\cal L}_H$
(but not a Hecke symmetry any more) satisfying relation $(3.1)$
at least for $T=S$ and $T\in\varepsilon(H_n(q))$. Since these arrows
generate the representation category of $S_\mu U(d)$ just as a
$C^*$--category, we deduce that equation $(3.1)$ holds on that
whole category.\medskip

We will see, as a consequence of Lemmas  5.1 and 5.4 together, that
for $n<d$ the
elements 
$S_{i_1,\dots, i_n}$ do {\sl not} belong to the maximal braided 
subcategory ${\cal L}_H^{\tilde{\varepsilon}}$.

Another consequence of Lemma 5.1 is that part e) allows one to construct
left inverses of $H$.
\medskip

\noindent{\bf 5.3 Corollary} {\sl Let $q>0$, and let, as before,  ${\cal
L}_H$ be the 
tensor $^*$--category of  Hilbert spaces with objects
tensor powers of the $d$--dimensional Hilbert space $H$. 
 Then, for $i_1<\dots<i_n$, the map 
$$\Phi^{i_1,\dots,i_n}:(H^{\otimes r}, H^{\otimes r})\to(H^{\otimes(r-1)},
H^{\otimes(r-1)}),$$
$$\Phi^{i_1,\dots,i_n}_r(T):=
\frac{1}{n!_{q}}S_{i_1,\dots,i_n}^*\otimes
1_{H^{\otimes(r-1)}}1_{H^{\otimes n-1}}\otimes 
TS_{i_1,\dots,i_n}\otimes 1_{H^{\otimes(r-1)}}$$
 is a positive left inverse of $H$.
If   $n=d$, $\Phi^{1,\dots,d}_2(g_q)=\lambda_{-d}1_H$.
If $n<d$, but $i_n=d$ and  $q\neq 1$, then $\Phi^{i_1,\dots,i_n}_2(g_q)$
is an invertible diagonal operator of $(H,H)$. In both cases these left inverses are faithful
on the range algebras of the Jimbo-Woronowicz Hecke symmetry
$\varepsilon: H_n(q)\to(H^{\otimes n}, H^{\otimes n})$.
}  
\medskip

\noindent{\bf Proof} Clearly  $\Phi^{i_1,\dots, i_n}_r$ is a positive
unital linear map with range 
 contained in the space of linear maps on $H^{\otimes(r-1)}$.
For $S\in(H^{\otimes r}, H^{\otimes r})$ and 
$T\in(H^{\otimes (r-1)}, H^{\otimes(r-1)})$,
$$\Phi^{i_1,\dots,i_n}_r (S1_H\otimes T)=\Phi^{i_1,\dots,i_n}_r(S)T,$$
so it is a left inverse of $H$. The Jimbo-Woronowicz Hecke $^*$--symmetry
is a braided symmetry
for ${\cal L}_H$ for which the range of $H_r(q)$ is a $^*$--subalgebra
of  $(H^{\otimes r},
 H^{\otimes r})^\varepsilon$, therefore whenever $\Phi^{i_1,\dots, i_n}_2(g_q)$ is 
invertible, $\Phi^{i_1,\dots, i_n}_r$ is faithul on that range, by Lemma
3.1.
\medskip

\noindent{\it Remark} Let $q>0$, and let $H(q)$ be endowed with its
standard involution. We have seen in section 3 that if we have a  Hecke
$^*$--symmetry $\varepsilon: H(q)\to{\cal T}$ in a tensor $C^*$--category with a
left inverse $\Phi$ for $\rho$ such that $\Phi(\varepsilon(g_1))$ is a nonzero scalar,
then the kernel of the $^*$--homomorphisms $\varepsilon: H_n(q)\to(\rho^n,\rho^n)$
are determined by that scalar. Therefore if $\Phi(\varepsilon(g_1))=\lambda_{-d}$,
$\varepsilon$ has dimension $d$.
\medskip

The following computations will serve 
to 
define a
certain tensor $C^*$--category with conjugates in the sense
of \cite{LR}.\medskip

\noindent{\bf 5.4 Lemma} {\sl Assume $q>0$. 
For $n\leq d$, and indices $i_1<\dots<i_n$ in ${1,\dots,d}$
consider the $n$--dimensional subspace $H_{i_1,\dots,i_n}$ of $H$
generated by $\psi_{i_1},\dots,\psi_{i_n}$. Then
$$S_{i_1,\dots,i_n}^*\otimes 1_H\circ 1_H\otimes
S_{i_1,\dots,i_n}=(n-1)!_{q}(-{\mu})^{n-1}1_{H_{i_1,\dots,i_{n}}},\eqno(5.1)$$
$$S_{i_1,\dots,i_{n}}^*\otimes 1_{H^{\otimes(n-1)}}\circ
1_{H^{\otimes(n-1)}}\otimes
{S_{i_1,\dots,i_n}}=(-{\mu})^{n-1}\varepsilon(A_{n-1})
1_{{H_{i_1,\dots,i_{n}}}^{\otimes 
(n-1)}}.\eqno(5.2)$$
In particular,
$$S^*\otimes 1_H\circ 1_H\otimes
S=(d-1)!_{q}(-{\mu})^{d-1}1_H,\eqno(5.3)$$
$$S^*\otimes 1_{H^{\otimes(d-1)}}\circ
1_{H^{\otimes(d-1)}}\otimes
{S}=(-{\mu})^{d-1}\varepsilon(A_{d-1}).\eqno(5.4)$$}
\medskip

\noindent{\bf Proof}  Set
$$\hat{\psi}_j=
\sum_{p(1)=j}(-\mu)^{i(p)}\psi_{i_{{p(2)}}}\otimes\dots\otimes
\psi_{i_{p(n)}}$$
and 
$$\tilde{\psi}_h=\sum_{p(d)=h}(-\mu)^{i(p)}\psi_{i_{{p(1)}}}
\otimes\dots\otimes\psi_{i_{{p(n-1)}}}$$
and write
$$S=\sum_j\psi_{i_j}\otimes\hat{\psi}_j=
\sum_h\tilde{\psi}_h\otimes\psi_{i_h}.$$
The left hand side of $(5.1)$ is the operator $T$ on $H$ with initial
and final support in $H_{i_1,\dots,i_n}$ such
that
$(\psi_{i_h},T\psi_{i_j})=(\hat{\psi}_j,\tilde{\psi}_h)$. Now,
 for $j\neq h$, $(\hat{\psi}_j,\tilde{\psi}_h)=0$, while 
$$(\hat{\psi}_j,\tilde{\psi}_j)=\sum_{p(1)=p'(n)=j}(-{\mu})^{i(p)}
(-\mu)^{i(p')}(\psi_{i_{p(2)}},\psi_{i_{p'(1)}})\dots
(\psi_{i_{p(n)}},\psi_{i_{p'(n-1)}}).$$
Thus for a fixed $p$ for which $p(1)=j$ we need to choose $p'$ so that
$p'(1)=p(2),\dots,p'(n-1)=p(n), p'(n)=j$. Since $i(p')=i(p)+n+1-2j$,
$$(\hat{\psi}_j,\tilde{\psi}_j)=
\sum_{p(1)=j}(-{\mu})^{i(p)}(-\mu)^{i(p)+n+1-2j}.$$
Writing each such $p$ in the form $p=(j-1j)\dots(23)(12)p''$ with
$p''(1)=1$, we see that $i(p)=i(p'')+j-1$, so $i(p')=i(p'')+n-j$
and 
$$(\hat{\psi}_j,\tilde{\psi}_j)=(n-1)!_{q}(-{\mu})^{n-1}.$$
 Taking into account the elements $\hat{\psi}_j$ and $\tilde{\psi}_h$,
the left hand side of 
$(5.2)$ is
$\sum_j\theta_{\hat{\psi}_j,\tilde{\psi}_j}$,
with $\theta_{\xi,\eta}(\zeta)=(\eta,\zeta)\xi.$
Defining, for each permutation $p$ for which $p(n)=j$, the permutation
$p'$ such that $p'(1)=j$, $p'(2)=p(1),\dots,p'(n)=p(n-1)$ shows that
$i(p)=i(p')+n+1-2j$, so $\tilde{\psi}_j=(-\mu)^{n+1-2j}\hat{\psi}_j$,
and we deduce that the left hand side of $(5.2)$ 
coincides
with
$\sum_j(-{\mu})^{n+1-2j}\theta_{\hat{\psi_j},\hat{\psi_j}}$.
Notice that
$(\hat{\psi}_j,\hat{\psi}_j)=(n-1)!_{q}q^{j-1}$.
Let $(j)$ denote the ordered $n-1$-tuple $1,\dots,n$ with the index 
$j$ suppressed. Then the relation with the elements $S_{(j)}$ 
 defined in  lemma $(5.1)$,
is 
$\hat{\psi}_j=(-\mu)^{j-1}S_{(j)}$,
therefore, in conclusion, $\sum_j\theta_{\hat{\psi}_j,\tilde{\psi}_j}=
(-{\mu})^{n-1}\sum{S}_{(j)}{S_{(j)}}^*=
(-{\mu})^{n-1}\varepsilon(A_{n-1})1_{{H_{i_1,\dots,i_n}}^{\otimes n-1}}.$
\medskip

{\it Remark} Notice that, for $n<d$, $q\neq 1$ and any choice of
indices $1\leq
i_1<\dots<i_n\leq d$, $1_H\otimes S_{i_1,\dots, i_n}$ {\sl is not}
a scalar multiple of $\varepsilon(g_1\dots g_n)S_{i_1\dots
i_n}\otimes1_H$.  
We show the claim.
Assume first that $(i_1,\dots, i_n)\neq (1,\dots, n)$. 
Then parts a) and e)
allow us to compute 
$$S_{i_1\dots i_n}^*\otimes 1_H\varepsilon(g_1\dots
g_n)S_{i_1\dots i_n}\otimes 1_H=q^{n-1}S_{i_1,\dots, i_n}^*\otimes 1_H 
1_{H^{\otimes
n-1}}\otimes g_qS_{i_1,\dots, i_n}\otimes 1_H$$ which must be a sum of
scalar
multiples of orthogonal
projections where,  for some $h>0$, a nonzero multiple of $1_{H^{(h)}}$
does appear.
On the other hand by Lemma 5.4 $S_{i_1,\dots, i_n}^*\otimes 1_H 1_H\otimes
S_{i_1,\dots, i_n}$ is just a scalar multiple of $1_{H_{i_1,\dots, i_n}}$.
In the case where $(i_1,\dots, i_n)=(1,\dots,n)$, arguing as in the proof 
of d) 5.1 shows that 
$$\varepsilon(g_1\dots g_n)S_{1,\dots, n}\otimes 1_{H^{(0)}}=(-\mu)^n
1_{H^{(0)}}\otimes S_{1,\dots, n}$$
which together with
$$\varepsilon(g_1\dots g_n)S_{1,\dots, n}\otimes 1_{H_{1,\dots
n}}=-(-\mu)^{n-1}
1_{H_{1,\dots, n}}\otimes S_{1,\dots, n}$$
and the fact that
$1_H=1_{H_{1,\dots, n}}+1_{H^{(0)}}$,
shows the claim
\medskip

Recall that  we defined the idempotents
$E_1=1$, $E_n=\frac{1}{n!_q}A_n$
of the Hecke algebra $H_\infty(q)$.
\medskip

\noindent{\bf 5.5 Theorem} {\sl Assume $q>0$.
Set, for $n\leq d$ and  indices $i_1<\dots< i_n$ in ${1,\dots,d}$,
$\overline{H}_{i_1,\dots,i_n}:=\varepsilon(E_{n-1})
{H_{i_1,\dots,i_n}}^{\otimes
(n-1)}$.
Then $$R_{i_1,\dots,i_n}:=
\frac{1}{\mu^{(n-1)/2}\sqrt{(n-1)!_q}}S_{i_1,\dots,i_n}\in
\overline{H}_{i_1,\dots,i_n}\otimes H_{i_1,\dots,i_n}$$
and
$$\overline{R}_{i_1,\dots,i_n}:=(-1)^{n-1}R_{i_1,\dots,i_n}=
\varepsilon'(n-1,1)R_{i_1,\dots,i_n}\in
H_{i_1,\dots,i_n}\otimes\overline{H}_{i_1,\dots,i_n}$$
satisfy the conjugate  equations:
$$\overline{R}_{i_1,\dots,i_n}^*\otimes 1_{H_{i_1,\dots,i_n}}\circ
1_{H_{i_1,\dots,i_n}}\otimes
R_{i_1,\dots,i_n}=1_{H_{i_1,\dots,i_n}},$$
$$R_{i_1,\dots,i_n}^*\otimes 1_{\overline{H}_{i_1,\dots,i_n}} \circ
1_{\overline{H}_{i_1,\dots,i_n}}\otimes
\overline{R}_{i_1,\dots,i_n}=1_{\overline{H}_{i_1,\dots,i_n}},$$
so $\overline{H}_{i_1,\dots,i_n}$ is a conjugate object for
$H_{i_1,\dots,i_n}$ in the smallest
tensor $^*$--subcategory of ${\cal L}_H$ containing 
$S_{i_1,\dots,i_n}$.}\medskip

\noindent{\bf Proof} We give a proof in the case $n=d$. In the general case, it suffices
to replace $H$ by the subspace $H_{i_1,\dots, i_n}$.
 We keep the notation of the proof of the two
previous 
lemmas. By
b) in Lemma 5.1, the range
of $\varepsilon(A_{d-1})$ is the linear span of $\{S_{(1)},\dots,
S_{(d)}\}$. Since each $S_{(j)}$ is a multiple of both
$\hat{\psi}_j$ and $\tilde{\psi}_j$,
${R}\in\overline{H}\otimes H$ and $\overline{R}\in H\otimes\overline{H}$. 
The proof of the  conjugate equations is immediate. 
\medskip

For future reference we present
a generalization of  relations $(5.3)$ and $(5.4)$ to
intermediate cases. 
\medskip

\noindent{\bf 5.6 Lemma} {\sl
For all
$k=1,\dots,d-1$,
$$S^*\otimes 1_{H^{\otimes k}}\circ 1_{H^{\otimes k}}
\otimes S=(d-k)!_qk!_{q}(-\mu)^{k(d-k)}\varepsilon(E_{k}).\eqno(5.5)
$$
Therefore the subobject of $H^{\otimes k}$ defined by $\varepsilon(E_k)$
has 
the subobject of $H^{\otimes(d-k)}$ defined by $\varepsilon(E_{d-k})$ as a 
conjugate 
defined by 
$$R_k=\frac{1}{\mu^{k(d-k)/2}\sqrt{k!_q}\sqrt{(d-k)!_q}} S,$$ 
$$\bar{R}_k=(-1)^{k(d-k)}R_k=\varepsilon'(d-k,k)R_k.$$
}
\medskip 

\noindent{\bf Proof} In order to simplify the following notation we
shall prove $(5.5)$ with $d-k$ in place of $k$.
Set
$$\tilde{\psi}_{i_1,\dots,i_k}=\sum_{p(d-k+1)=i_1,\dots,p(d)=i_k}
(-\mu)^{i(p)}\psi_{p(1)}\otimes\dots\otimes\psi_{p(d-k)}.$$
Then $$S=\sum_{i_1<\dots<i_k}\sum_{p\in{\Bbb
P}_k}
\tilde{\psi}_{i_{p(1)},\dots,i_{p(k)}}\otimes
\psi_{i_{p(1)}}\otimes\dots\otimes\psi_{i_{p(k)}}.$$
Now, if $i_1<\dots<i_k$,
$$\tilde{\psi}_{i_{p(1)},\dots,i_{p(k)}}=
(-\mu)^{i(p)}\tilde{\psi}_{i_1,\dots,i_k},$$
so
$$S=\sum_{i_1<\dots<i_k}\tilde{\psi}_{i_1,\dots,i_k}\otimes 
S_{i_1,\dots,i_k},$$
as $i(p)=i(p^{-1})$.
Similarly,
${S}=\sum_{i_1<\dots<i_k}
{S}_{i_1,\dots,i_k}\otimes
\overline{\psi}_{i_1,\dots,i_k}$ with
$$\overline{\psi}_{i_1,\dots,i_k}=\sum_{p(1)=i_1,\dots,p(k)=i_k}
(-{\mu})^{i(p)}\psi_{p(k+1)}\otimes\dots\otimes\psi_{p(d)}.$$
Thus the left hand side of $(5.5)$, with, recall,  $k$ and $d-k$ replaced,
equals
$$\sum_{i_1<\dots<i_k}
({S}_{i_1,\dots,i_k},S_{i_1,\dots,i_k})
\theta_{\tilde{\psi}_{i_1,\dots,i_k},\overline{\psi}_{i_1,\dots,i_k}},$$
where
$({S}_{i_1,\dots,i_k},S_{i_1,\dots,i_k})=k!_{q}$.
Now 
$\sum_{i_1<\dots<i_k}
\theta_{\tilde{\psi}_{i_1,\dots,i_k},\overline{\psi}_{i_1,\dots,i_k}}$
annihilates all simple tensors
$\psi_{j_1}\otimes\dots\otimes\psi_{j_{d-k}}$
with $j_h=j_k$ for some $h\neq k$, and applied on a  simple tensor
of the form $\psi_{j_{p(1)}}\otimes\dots\otimes\psi_{j_{p(d-k)}}$, with
$j_1<\dots<j_{d-k}$ and $p\in{\Bbb P}_{d-k}$, yields
$(-\mu)^{i(p)+N}  \tilde{\psi}_{i_1,\dots,i_k}$, with
$(i_1,\dots,i_k)$ obtained from $(1,\dots,d)$ removing the indices
$j_1,\dots, j_{d-k}$, and with $N$ the number of pairs $(i_h, j_k)$ for
which $i_h>j_k$. On the other hand it is not difficult to check that
$$\tilde{\psi}_{i_1,\dots,
i_k}=S_{j_{1},\dots,j_{d-k}}(-\mu)^{k(d-k)-N},$$ 
so the image of the left hand side of $(5.5)$ on that simple 
tensor is 
$$(-\mu)^{i(p)}k!_q(-\mu)^{k(d-k)}S_{j_1,\dots,j_{d-k}},$$
the same as the image under $(-\mu)^{k(d-k)}k!_q\varepsilon(A_{d-k}).$ 
\medskip

\noindent{\it Remark} 
Choosing $k=2$ shows that $\varepsilon(E_2)$, and therefore
 the Hecke symmetry $\varepsilon$, is determined uniquely 
by the intertwiner $S$ in the tensor $C^*$--category ${\cal L}_H$.
\medskip

Let us consider  the smallest
tensor $^*$--subcategory of ${\cal L}_H$ with 
 subobjects, containing  $S$. This is the representation category
$\text{Rep}(S_\mu U(d))$
of the Woronowicz compact quantum group $S_\mu U(d)$.
It is a tensor $C^*$--category of Hilbert spaces with conjugates
and  with a Hecke
 symmetry whose dimension equals the dimension of $H$. 
\bigskip

\end{section}

\begin{section} {An abstract characterization of $\text{Rep}(S_\mu U(d))$}

In the previous section we have shown the following relations between the
fundamental intertwiner $S$ of the representation category of the Woronowicz quantum group
$S_\mu U(d)$ and the model Woronowicz--Jimbo q-Hecke $^*$--symmetry, with 
$q=\mu^2$.
$$S^*S=d!_q,\eqno(6.1)$$
$$S^*\otimes 1_H\circ 1_H\otimes S=(d-1)!_q(-\mu)^{d-1},\eqno(6.2)$$
$$SS^*=\varepsilon(A_d),\eqno(6.3)$$
$$\varepsilon(g_1\dots g_d)S\otimes 1_H=-(-\mu)^{d-1}1_H\otimes
S.\eqno(6.4)$$
 The second and third  relations above serve to define
conjugates
 in the tensor $C^*$--category 
$\text{Rep}(S_\mu U(d))$, while the last one asserts that 
that category is 
a braided tensor category with  respect to a renormalization of the 
Woronowicz-Jimbo Hecke
symmetry.

Let  ${\cal T}$  denote a tensor $C^*$--category
with tensor unit $\iota$ 
such that $(\iota,\iota)={\Bbb C}$
 and objects $\{\iota,\rho,\rho^2,\dots\}$, the
tensor powers of a fixed  object $\rho$. 
We shall also assume the existence of a Hecke $^*$--symmetry
$\varepsilon: H_n(q)\to(\rho^n,\rho^n)$ for some $q>0$. 

In this framework we introduce the notion of {\sl special object}
which generalizes the corresponding notion in permutation 
symmetric tensor $C^*$--categories
given by Doplicher and Roberts in \cite{DR}.\medskip

\noindent{\bf 6.1 Definition} The  object
$\rho$ of ${\cal T}$ is called {\sl special with dimension}
$d$ if there is an intertwiner $R\in(\iota,\rho^d)$
such that properties $(6.1)$--$(6.4)$ hold for $R$ and $\rho$ in place
of $S$ and $H$ respectively and $\mu$ the positive square root of $q$.
 $R$ will be referred to as a {\sl special
intertwiner}.
\medskip

The aim of this section is to show that the above properties are  sufficient to
construct a embedding of tensor $C^*$--categories
from ${\text Rep}(S_\mu U(d))$ to ${\cal T}$.
\medskip

\noindent{\bf 6.2 Theorem} {\sl Let  $\rho$ be a special 
object of ${\cal T}$
with dimension  $d\geq2$.
 Then there is a unique faithful
tensor $^*$--functor $\text{Rep}(S_\mu U(d))\to {\cal T}$
taking the fundamental representation  of $S_\mu U(d)$ to $\rho$
and $S$ to a special intertwiner $R$. } \medskip

\noindent{\bf Proof} Notice that for $q=1$, this theorem reduces to
Theorems 4.1 and 4.4 in \cite{DRduals}. Here we argue in the same way. 
We first show that our assumptions imply that $\varepsilon$
has dimension $d$. By Theorem 3.3 we 
need to show that $\rho$ has a positive left inverse $\Phi$ such that
$\Phi(\varepsilon(g_1))=\lambda_{-d}$. In fact,
$$\Phi(T):=\frac{1}{d!_q}R^*\otimes 1_\rho\circ 1_{\rho^{d-1}}\otimes
T\circ
R\otimes 1_\rho$$ is a
positive, unital left
inverse of $\rho$ such that, 
since  by $(6.1)$, $(6.3)$ and Lemma 2.1 a), $\varepsilon(g_j)R=qR$,
$j=1,\dots,d-1$,
$$\Phi(\varepsilon(g_1))=\frac{1}{d!_q}R^*\otimes
1_\rho(\varepsilon(g_1\dots
g_{d-1}))^{-1}\varepsilon(g_1\dots g_{d})R\otimes 1_\rho=$$
$$
\frac{1}{q^{d-1}d!_q}
R^*\otimes 1_\rho\varepsilon(g_1\dots g_d)R\otimes
1_\rho=-\frac{1}{(-\mu)^{d-1}d!_q}
R^*\otimes 1_\rho 1_\rho\otimes
R=$$
$$-\frac{1}{1+q+\dots+q^{d-1}}=\lambda_{-d}.$$
Therefore the Hecke symmetry $H(q)\to{\cal T}$ factors over
the model symmetry, and this defines a faithful tensor
$^*$--functor from homogeneous  subcategory (the subcategory with
the same objects and with intertwiners  with  the same domain and
range) 
of
$\text{Rep}(S_\mu U(d))$
to ${\cal T}$. Now any intertwiner space 
$(H^{\otimes r}, H^{\otimes s})$
in 
$\text{Rep}(S_\mu U(d))$
is nonzero if and only if $s-r$ is an integer multiple of $d$, and, for
$s=r+kd$, $k>0$, one can write
$T\in(H^{\otimes r}, H^{\otimes s})$ uniquely in the form
$T=T_0\circ S\otimes\dots\otimes S$ with
$T_0\in
(H^{\otimes s}, H^{\otimes s})$ satisfying
$T_0=T_0\circ SS^*\otimes\dots\otimes SS^*.$
One can thus uniquely extend that functor from the homogeneous
subcategory  to a
faithful tensor
$^*$--functor on the whole of $\text{Rep}(S_\mu U(d))$, taking $S$ to $R$. 
\bigskip

\noindent{\sl Acknowledgements} The author thanks S. Doplicher and J.E.
Roberts for discussions and for their the constant interest  during the
preparation of this work.

\end{section}

 \end{document}